\documentclass[12pt]{article}
\usepackage[utf8]{inputenc}
\usepackage[T1]{fontenc}
\usepackage{amsmath,amssymb,amsfonts,amsthm,enumerate}
\usepackage{mathrsfs}
\usepackage{graphicx}
\usepackage{multicol,multirow}
\usepackage{rotating}
\usepackage{appendix}
\usepackage{textcomp}

\newtheorem{theorem}{Theorem}[section]
\newtheorem{lemma}[theorem]{Lemma}
\newtheorem{defin}[theorem]{Definition}
\newtheorem{remark}[theorem]{Remark}

\numberwithin{equation}{section}
\newtheorem{proposition}[theorem]{Proposition}

\newtheorem{corollary}[theorem]{Corollary}

\def\rr{{\mathbb R}}
\def\cc{{\mathbb C}}
\def\qq{{\mathbb Q}}
\def\zz{{\mathbb Z}}
\def\ee{{\cal E}}

\def\rkk{\rr ^k}
\def\rk1{\rr ^{k+1}}
\def\su{\subset}

\def\al{\alpha}
\def\be{\beta}
\def\ga{\gamma}
\def\Ga{\Gamma}

\def\ep{\varepsilon}
\def\si{\sigma}
\def\Si{\Sigma}
\def\la{\lambda}

\def\Om{\Omega}
\def\cd{\cdot}
\def\stb{,\ldots ,}
\def\emp{\emptyset}
\def\ol{\overline}

\def\ii{{\cal I}}

\def\ff{{\cal F}}
\def\sik{\rr ^2}
\def\sumjn{\sum_{j=1}^n}
\def\sumjs{\sum_{j=1}^s}
\def\ran{\rangle}
\def\lan{\langle}
\def\proof{\noindent {\bf Proof.} }
\def\qued{\hfill $\square$}

\begin{document}

\title{Solutions to the discrete
Pompeiu problem and to the finite Steinhaus tiling problem}

\author{Gergely Kiss and Mikl\'os Laczkovich}



\maketitle

\begin{abstract}
Let $K$ be a nonempty finite subset of the Euclidean space $\rkk$
$(k\ge 2)$. We prove that if a function $f\colon \rkk \to \cc$ is such that the
sum of $f$ on every congruent copy of $K$ is zero, then $f$ vanishes everywhere.
In fact, a stronger, weighted version is proved. As a corollary we find that
every finite subset of $\rkk$ having at least two elements is a Jackson
set; that is, no subset of $\rkk$ intersects every congruent copy of $K$ in
exactly one point.
\end{abstract}

\section{Introduction and main results}
A compact subset $K$ of the plane having positive Lebesgue measure is said
to have the {\it Pompeiu property} if the following condition is
satisfied: if $f\colon \sik \to \cc$ is a continuous function such that
$\int_{\si (K)} f\, d\la _2 =0$ for every isometry $\si$ of the plane, then
$f\equiv 0$. It is known that the disc does not have the Pompeiu property, while
all polygons have. The classical Pompeiu problem asks if a connected compact
set with a smooth boundary that does not have the Pompeiu property is
necessarily a disc. As for the history of the problem, see \cite{KS}, \cite{MR},
\cite{R}, \cite{Za}.

Replacing the Lebesgue measure $\la _2$ by the counting measure and the
isometry group by an arbitrary family $\ff$ of functions mapping a nonempty
set $X$ into itself, we obtain the following notion.

\begin{defin}
{\rm
Let $\ff$ be a family of functions mapping a set $X$ into itself, and let
$K$ be a nonempty finite subset of $X$. We say that} $K$ has the
discrete Pompeiu property with respect to the family $\ff$ if the following
condition is satisfied: whenever $f\colon X\to \cc$ is such that
\begin{equation}\label{e8}
\sum_{x\in K} f(\phi (x))=0
\end{equation}
for every $\phi \in \ff$, then $f\equiv 0$.
\end{defin}
We also define a stronger property as follows.
\begin{defin}
We say that an $n$-tuple $(a_1 \stb a_n )$ of (not necessarily distinct)
elements of $X$ has the weighted discrete Pompeiu property with respect to
the family $\ff$ if the following condition is satisfied: whenever
$c_1 \stb c_n$ are complex numbers with $\sumjn c_j \ne 0$ and $f\colon X\to
\cc$ is such that $\sumjn c_j \cd f(\phi (a_j ))=0$ for every $\phi \in \ff$,
then $f\equiv 0$.
\end{defin}

Note that the condition $\sumjn c_j \ne 0$ is necessary: if $\sumjn c_j =0$,
then every constant function satisfies the condition.

Previous results in this topic focused on the following families $\ff$:
\begin{enumerate}[$\bullet$]
\item translations in groups;
\item similarities of $\rkk$; that is, maps $\phi \colon \rkk \to \rkk$
satisfying $|\phi (x)-\phi (y)|=q\cd |x-y|$ for every $x,y\in \rkk$ with a
suitable $q>0$;
\item direct similarities of $\rkk$; that is, similarities preserving
orientation;
\item isometries (or congruences) of $\rkk$; that is, distance preserving
bijection of $\rkk$ onto itself;
\item rigid motions of $\rkk$; that is, isometries preserving orientation.
\end{enumerate}

Note that in $\rkk$ the families of translations, similarities, direct
similarities and rigid motions are listed in descending order. Therefore, if
a subset $K\su \rkk$ has the (weighted) discrete Pompeiu property with respect
to any of these families, then it has the same property with respect to the
previous families as well. As for the history of the subject up to 2018, we
quote \cite[p.~300]{KLV}.

``Apparently, the first results concerning the discrete
Pompeiu property appeared in \cite{Z}, where the author considers the Pompeiu
problem for finite subsets of $\zz ^n$ with respect to  translations. The
interest in the topic revived shortly after the 70th William Lowell Putnam
Mathematical  Competition (2009), where the following problem was posed: Let
$f$ be a real-valued function on the plane such that for every square $ABCD$
in the plane, $f(A)+f(B)+f(C)+f(D)=0.$ Does it follow that $f\equiv 0$? This is
nothing but asking whether the set of vertices of a square has the discrete
Pompeiu property with respect to the similarities of the plane. This problem
motivated the paper \cite{GD} by C. De Groote and M. Duerinckx. They proved
that every finite and nonempty subset of $\sik$ has the discrete Pompeiu
property with respect to direct similarities. Another generalization of the
Putnam problem appeared in \cite{KKS}, where it is proved that the set of
vertices of the unit square has the discrete Pompeiu property with respect to
the group of isometries. Recently, M. J. Puls \cite{Puls} considered the
discrete Pompeiu problem in groups.''

The case of translations in groups is also treated in \cite{KLV}, \cite{KMS}
and \cite{LP}. According to \cite[Proposition 1.1]{KLV}, no finite subset
having at least two elements of a torsion free Abelian group (in particular, of
$\rkk$) has the discrete Pompeiu property with respect to translations.

The weighted version of the result of C. De Groote and M. Duerinckx was proved
in \cite[Theorem 3.3]{KLV}.

The notion of similarity can be defined in every Abelian group as follows.
We say that the map $\phi \colon G\to G$ is a {\it simple similarity} of the
Abelian group $G$ if there is an element $a\in G$
and there is a positive integer $k$ such that $\phi (x)=a+k\cd x$ for every
$x\in G$. The following statement generalizes the results of \cite{GD} and
\cite{KLV} cited above.
\begin{proposition}\label{p1}
In every Abelian group G, every $n$-tuple of points of $G$ has the weighted
discrete Pompeiu property with respect to the family of simple similarities.
\end{proposition}

In the sequel we consider the case when $X=\rkk$ and $\ff =G_k$ is the family
of rigid motions of $\rkk$. As we mentioned above, in this context the first
relevant result appeared in \cite{KKS}, stating that the set of vertices of
the unit square has the discrete Pompeiu property with respect to $G_2$. Later
the discrete Pompeiu property of all parallelograms and
some other special four-element subsets of the plane was established in
\cite{KLV}. Our main result is the following.

\begin{theorem}\label{t1}
For every $k\ge 2$ and $a_1 \stb a_n \in \rkk$, the $n$-tuple $(a_1 \stb a_n )$
has the weighted discrete Pompeiu property with respect to the group $G_k$ of
rigid motions of $\rkk$. In particular, every nonempty finite subset of
$\rkk$ has the discrete Pompeiu property with respect to $G_k$.
\end{theorem}

In fact we prove the following stronger result.
\begin{theorem}\label{t1a}
Let $k\ge 2$ and $a_1 \stb a_n \in \rkk$ be given. Then there is a
countable additive subgroup $E$ of $\rkk$ with the following property: whenever
$c_1 \stb c_n \in \cc$, $\sumjn c_j \ne 0$, and $f\colon E\to \cc$ is such that
$\sumjn c_j \cd f(\phi (a_j ))=0$ for every $\phi \in G_k$ satisfying
$\phi (a_j )\in E$ for every $j=1\stb n$, then $f\equiv 0$.
\end{theorem}
We can deduce Theorem \ref{t1} from Theorem \ref{t1a} as follows.

\proof
Let $c_1 \stb c_n$ be complex numbers with $\sumjn c_j \ne 0$, and suppose that
$f\colon \rkk\to \cc$ satisfies $\sumjn c_j \cd f(\phi (a_j ))=0$ for every
$\phi \in G_k$. Let $E$ be as in Theorem \ref{t1a}, and let $z_0 \in \rkk$ be
arbitrary. Then the function $g(x)=f(x+z_0 )$ $(x\in E)$ satisfies the
conditions of Theorem \ref{t1a}. Therefore, we have $g\equiv 0$, and thus
$f(z_0 )=0$.
\qued
\begin{remark}\label{r1}
{\rm
For $k=1$ the statement of Theorem \ref{t1} is false: if $n\ge 2$,
$K=\{ 1\stb n\}$ and $f(x)=e^{2\pi ix/n}$, then \eqref{e8} holds for every
$\phi \in G_1$. More generally, if $K=\{ a_1 \stb a_n \}$, $\la$ is a root of
the entire function $\sumjn e^{a_j z}$ and $f(x)=e^{\la x}$, then
\eqref{e8} holds for every translation $\phi$. If $K$ is symmetric; that is,
if $K=-K$, then \eqref{e8} holds for every isometry $\phi$ of $\rr$.}
\end{remark}

Let $a_1 \stb a_n$ and $z_0$ be given points in $\rkk$, and let $c_1 \stb c_n$
be complex numbers with $\sumjn c_j \ne 0$. Let $\Si$ denote the system
of linear equations
$$\sumjn c_j \cd x_{\phi (a_j )} =0 \quad (\phi \in G_k ),\qquad
x_{z_0}=1,$$
where $x_z$ is an unknown for every $z\in \rkk$.
By Theorem \ref{t1}, $\Si$ has no solution. Now, it is easy to see that
if a system of linear equations has no solution, then there is a finite
subsystem of $\Si$ that has no solution either\footnote{For the sake of
completeness we provide the simple proof in Section 6.}.
Therefore, we obtain the following.
\begin{corollary}\label{c1}
Suppose $k\ge 2$, $a_1 \stb a_n ,z_0 \in \rkk$ and $c_1 \stb c_n \in \cc$
are given such that $\sumjn c_j \ne 0$. Then there is a finite set $H\su \rkk$
containing $z_0$ with the following property: whenever $f\colon H\to \cc$ is
such that  $\sumjn c_j \cd f(\phi (a_j ))=0$ for every $\phi \in G_k$
satisfying $\phi (a_j )\in H$ for every $j=1\stb n$, then $f(z_0 )=0$.
\end{corollary}

\subsection*{The structure of the paper} Before turning to the proofs,
in the next section we discuss some consequences of Theorem \ref{t1} in
the Euclidean Ramsey theory and in the topic of Steinhaus sets. (For the
relevant notions see therein.) In Section 3 we explain an important ingredient
of the proof of Theorem \ref{t1a}, namely harmonic analysis
on discrete groups. In addition, we prove Proposition \ref{p1} using this
method. In Sections 4 and 5 we prove Theorem \ref{t1a}
for $k=2$ and for $k>2$, respectively. Finally, Section 6 gives a proof of our
remark on infinite systems of linear equations.

\section{Applications to coloring problems and to the finite Steinhaus
tiling problem}
Theorem \ref{t1} has the following obvious consequence.
\begin{corollary}\label{c3}
If $k\ge 2$, $K\su \rkk$ has $n$ elements, $d\mid n$, $d>1$ and $\rkk$ is
colored with $d$ colors, then there is a congruent copy of $K$ containing more
than $n/d$ points of the same color.
\end{corollary}

\proof
If the statement of the corollary is not true, then there is a partition
$\rkk =A_1 \cup \ldots \cup A_d$ such that every congruent copy of $K$
intersects each of the sets $A_1 \stb A_d$
in exactly $n/d$ points. Let $b_1 \stb b_d$ be nonzero complex numbers with
$\sum_{j=1} ^d b_j =0$. If we define $f(x)=b_j$ for every $x\in A_j$
$(j=1\stb d)$, then $\sum_{x\in K} f(\phi (x))=\sum_{j=1} ^d (n/d)\cd b_j =0$
for every $\phi \in G_k$, contradicting Theorem \ref{t1}.
\qued

In the case of $n=4$, $d=2$ we obtain the following.
\begin{corollary}\label{c4}
If $k\ge 2$, $K\su \rkk$, $|K|=4$ and if $\rkk$ is colored with two colors, then there is a
congruent copy of $K$ containing at least three points of the same color.
\end{corollary}
If $k=2$ and $K$ is a rectangle, then
we obtain the following: {\it For every right triangle $T$ and for every
coloring the plane with two colors, there is always a monochromatic
triangle congruent to $T$.} This is L.E. Shader's theorem \cite{S}.

The special case $d=n$ in Corollary \ref{c3} is closely connected to
the general Steinhaus problem. Let $E\su \rkk$ be given. We say that $S\su \rkk$
is a {\it Steinhaus set for $E$}, if $S$ intersects every congruent copy
of $E$ in exactly one point. The original question of Hugo Steinhaus,
posed in the 1950s, was the following: is there a Steinhaus set in the
plane for $\zz ^2$? This question was answered in the affirmative by S.
Jackson and R.D. Mauldin in 2002 \cite{JM1} (see also \cite{JM}). Analogous
results were obtained by P. Komj\'ath \cite{K1}, \cite{K2} and J.H. Schmerl
\cite{Sch} for $\zz$, $\qq$ and $\qq ^n$.

These results motivated S. Jackson to ask {\it if there are sets $K,S\su \sik$
such that the cardinality of $K$ is finite but at least two, and $S$ is a
Steinhaus set for $K$.}
\begin{defin}
{\rm A finite set $K\su \rkk$ is called a} Jackson set {\rm if there is no
Steinhaus set for $K$.}
\end{defin}
(The notion of Jackson set was introduced in \cite{GMW}). It is clear that
singletons are not
Jackson sets (as $S=\rkk$ works), and it is easy to see that if $k\ge 2$,
then all $2$-element subsets of $\rkk$ are Jackson sets. It is known that every
subset of $\sik$ of cardinality $3,4,5$ or $7$ is a Jackson set (see
\cite{HJL}). It is also known that for every finite set $K\su \sik$ having at
least two elements there are no measurable Steinhaus sets for $K$; see
\cite{KP}.

\begin{proposition}\label{p2}
Let $K\su \rkk$ be a finite set of cardinality at least two. If $K$ has the
discrete Pompeiu property with respect to $G_k$, then $K$ is a Jackson set.
\end{proposition}

\proof
We apply the argument of \cite[Proposition 1.3]{GMW}. Suppose that
$K\su \rkk$, $|K|\ge 2$, and that $S\su \rkk$ is a Steinhaus set for $K$.
Then the sets $S-a$ $(a\in K)$ are pairwise disjoint. Indeed, if
$c\in (S-a)\cap (S-b)$, where $a,b\in K$ and $a\ne b$, then $c+a, c+b\in S$.
In this case, however,
$|S\cap \si (K)|\ge 2$ for the translation $\si (x)=x+c$, which is impossible.
We have $\bigcup_{a\in K} (S-a)=\rkk$. Indeed, if $x\in \rkk$ is arbitrary and
$S\cap (K+x)=\{ s\}$, then $s=a+x$, where $a\in K$, and thus $x=s-a\in S-a$.
Therefore, the sets $S-a$ $(a\in K)$ constitute a partition of $\rkk$ such
that every congruent copy of $K$ intersects each of the sets $S-a$ in exactly
one point. As we saw in the proof of Corollary \ref{c3}, this contradicts
the discrete Pompeiu property of the set $K$.
\qued

By Theorem \ref{t1} we obtain the following:
\begin{corollary}\label{c5}
Every finite subset of $\rkk$ $(k\ge 2)$ having at least two elements is a
Jackson set.
\end{corollary}

\begin{remark}\label{r4}
{\rm
For $k=1$ the statement of the corollary is false: if $K=\{ 1\stb n\}$, then
$S=\bigcup_{t\in \zz} ([0,1)+n\cd t)$ intersects every congruent copy of $K$
in exactly one point, so $K$ is not a Jackson set. For more on Jackson sets
in $\rr$, see \cite{GMW}.}
\end{remark}

\begin{remark}
{\rm
Note that the definition of Jackson set uses isometries and not just rigid
motions, while Theorem \ref{t1} is about the Pompeiu property with respect to
the family of rigid motions. Therefore, Corollaries \ref{c3}, \ref{c4},
\ref{c5} remain true if we replace congruent copies of $K$ by images of $K$
under rigid motions, and, in the definition of Jackson sets
we replace isometries by rigid motions.}
\end{remark}

\begin{remark}
{\rm
Let $K\su \rkk$ be given, and let $m$ be a positive integer. We say that the
set $S\su \rkk$ is an {\it $m$-Steinhaus set for $K$} if every congruent copy
of $S$ intersects $K$ in exactly $m$ points. The finite set $K$ is called an
{\it $m$-Jackson set}, if there is no $m$-Steinhaus set for $K$.
Obviously, the sets of cardinality $<m$ are $m$-Jackson sets, and if
$|K|=m$, then $K$ is not an $m$-Jackson set, as $S=\rkk$ is an
$m$-Steinhaus set for $K$. The following
generalization of Corollary \ref{c5} can be obtained by a similar argument.
\begin{corollary}\label{c6}
Every finite subset of $\rkk$ $(k\ge 2)$ having more than $m$ elements is an
$m$-Jackson set.
\end{corollary}

\proof
Suppose $S$ is an $m$-Steinhaus set for $K$. Then the sets
$S-a$ $(a\in K)$ constitute an $m$-cover of $\rkk$. Indeed, if $x\in \rkk$,
then $x\in S-a$ $(a\in K)\iff x+a \in S$. Since $|S\cap (K+x)|=m$, it follows
that every point of $\rkk$ is contained in exactly $m$ of the sets $S-a$
$(a\in K)$.

Let $|K|=n>m$, and put $f(x)=1/m$ if $x\in S$, and $f(x)=-1/(n-m)$ if
$x\notin S$. Then $f$ is nowhere zero, but $\sum_{x\in K} f(\si (x))=0$ for
every $\si \in G_k$. By Theorem \ref{t1}, this is impossible.
\qued
}
\end{remark}

\section{Harmonic analysis on discrete Abelian groups}
Similarly to the classical Pompeiu problem, the main tool in proving Theorem
\ref{t1} is harmonic analysis. Since our objects are finite, we need
harmonic analysis on discrete groups. Let $G$ be an Abelian group equipped with
the discrete topology. We denote by $C(G)$ the set of all maps from $G$ into
$\cc$ equipped with the topology of pointwise convergence. More precisely,
a set $U\su C(G)$ is open if, for every $f\in U$ there is a finite set $F\su G$
and an $\ep >0$ such that, if $g\in C(G)$ and if $|g(x)-f(x)|<\ep$ for
every $x\in F$, then $g\in U$. (In fact, this is the same as the product
topology of $\cc ^G$.) A nonzero function $m\in C(G)$ is called an
{\it exponential}, if $m$ is multiplicative; that is, if $m(x+y)=m(x)\cd m(y)$
for every $x,y\in G$. Note that if $m$ is an exponential, then $m(0)=1$,
and $m(x)\ne 0$ and $m(-x)=1/m(x)$ for every $x\in G$.

By a {\it variety} we mean a translation invariant closed linear subspace of
$C(G)$. We say that {\it harmonic analysis} holds on $G$ if every nonzero
variety contains an exponential. The following result will be needed in the
proof of Theorem \ref{t1a}.

\begin{theorem}\label{tha} {\rm \cite[Theorem 1]{LS}}
Harmonic analysis holds on a discrete Abelian group $G$ if and only if the
torsion free rank of $G$ is less than continuum.
\end{theorem}
According to this theorem, harmonic analysis does not hold on the additive
group of $\rkk$. On the other hand, it holds on every countable Abelian
group, and so we have to work on suitable countable
subgroups of $\rkk$. The next proof of Proposition \ref{p1} is hardly more
than an application of this fact.

\subsection*{Proof of Proposition \ref{p1}} Let $G$ be an Abelian group, and let
$a_1 \stb a_n \in G$ and $c_1 \stb c_n \in \cc$ be given such that $\sumjn
c_j \ne 0$. Let $f\colon G\to \cc$ be such that $\sumjn c_j \cd f(b+k\cd a_j )
=0$ for every $b\in G$ and $k=1,2, \ldots$. We have to prove that $f\equiv 0$.
Suppose this is not true, and let $x\in G$ be such that $f(x)\ne 0$.
Let $H$ denote the subgroup of $G$ generated by $x$ and $a_1 \stb a_n$,
and let $V$ denote the set of all functions $g\colon H \to \cc$ such that
$\sumjn c _j \cd g(b+k\cd a_j ) =0$ for every $b\in H$ and $k=1,2,\ldots $. It
is clear that $V$ is a linear space over $\cc$, and that $V$ is invariant
under translations by elements of $H$. It is also easy to see that $V$ is
closed in the set $\cc ^H$ equipped with the product topology. This means
that $V$ is a variety on the discrete, countable additive group $H$.

Since $f|_H \in V$, we have $V \ne \{ 0\}$. Then, by Theorem \ref{tha},
$V$ contains an exponential; that is, a function $m\colon H\to \cc$ such that
$m\ne 0$ and $m(x+y)=m(x)\cd m(y)$ for every $x,y\in H$.
Since $m\in V$, we have
\begin{equation}\label{e11}
\sumjn c_j \cd  m(a_j )^k =\sumjn c_j \cd  m(k\cd a_j )=0 \qquad
(k=1,2,\ldots ).
\end{equation}
Permuting the elements $a_1 \stb a_n$ if necessary, we may assume that
there is an $1\le s\le n$ such that $m(a_1 )\stb m(a_s )$ are distinct
nonzero complex numbers, and
for every $s<j\le n$, $m(a_j )$ equals one of $m(a_1 )\stb m(a_s )$.
Put $d_j =\sum \{ c_\nu \colon m(a_\nu )=m(a_j )\}$
for every $j=1\stb s$; then $\sumjs d_j =\sumjn c_j \ne 0$.

By \eqref{e11} we have $\sumjs d_j \cd  m(a_j )^k =0$ for every $k=1,2, \ldots$,
and thus
\begin{equation}\label{e12}
\sumjs d_j \cd m(a_j ) \cd m(a_j )^k =0 \qquad (k=0\stb s-1).
\end{equation}
Now, \eqref{e12} is a system of linear equations with
unknowns $d_j \cd m(a_j )$ $(j=1\stb s)$. The determinant of this system,
being a Vandermonde determinant, is nonzero. Therefore, we have
$d_1 =\ldots =d_s =0$, which is impossible. This contradiction proves the
statement. \hfill $\square$

\section{Proof of Theorem \ref{t1a} for $k=2$}
Let $(a_1 \stb a_n )$ be a fixed $n$-tuple of elements of $\sik$. We identify
$\sik$ with $\cc$ (the field of complex numbers), and denote by $S^1$ the unit
circle $\{ u\in \cc \colon |u|=1\}$. Let $E$ denote the subfield of $\cc$
generated by $a_1 \stb a_n$ and
the set $S^1_a =\{ z\in S^1 \colon z$ is algebraic$\}$. Then $E$ is a countable
subfield of $\cc$ containing $S^1_a \cup \{ a_1 \stb a_n \}$.

We show that the set $E$ satisfies the condition of Theorem \ref{t1a} (in
the case of $k=2$). More precisely, we prove that {\it if $c_1 \stb c_n \in
\cc$, $\sumjn c_j \ne 0$ and $f\colon E \to \cc$ is such that $\sumjn
c _j \cd f(x+a_j y) =0$ for every $x\in E$ and $y\in S^1_a$, then $f\equiv 0$.}

\subsection*{The structure of the proof} Suppose $f$
satisfies the condition, but $f\not\equiv 0$. Applying Theorem \ref{tha} on
the countable additive group $E$, we find an exponential
$m\colon E\to \cc$ such that $\sumjn c_j \cd m(u\cd a_j )=0$ for every
$u\in S^1_a$.  Then we apply this equation with many $u$ having rational
coordinates, and obtain, by applying a theorem of J.-H. Evertse,
H.P. Schlickewei and W.M. Schmidt on the number of solutions of linear
equations, that there is an integer $d$ such that for every algebraic $u$
with $|u|=1$ there are indices $j_1 \ne j_2$ such that $m(u\cd (a_{j_2}-a_{j_1}) /d)$ and $m(u\cd i\cd (a_{j_2}-a_{j_1}) /d)$ are roots of unity of bounded degree
(Lemma \ref{l1}). In the final step we show that this contradicts the
fact that $m(x_1 )\cdots m(x_s )=1$ whenever $x_1 +\ldots +x_s =0$.

\medskip
Now we turn to the details.
Fix $c_1 \stb c_n \in \cc$ with $\sumjn c_j \ne 0$. Clearly, we may assume that
$c _j \ne 0$ for every $j=1\stb n$. Let $\Om$ denote the set of all functions
$f\colon E \to \cc$ such that
$\sumjn c _j \cd  f(x+a_j y) =0$ for every $x\in E$ and $y\in S^1_a$. It
is clear that $\Om$ is a linear space over $\cc$, and that $\Om$ is invariant
under translations by elements of $E$. It is also easy to see that $\Om$ is
closed in the set $\cc ^E$ equipped with the product topology. This means
that $\Om$ is a variety on the discrete additive group $E$.

Suppose that the statement of the theorem is false; that is, $\Om \ne \{ 0\}$.
Clearly, this implies $n\ge 2$. Then, by Theorem \ref{tha}, $\Om$ contains an
exponential; that is, a function $m\colon E\to \cc$ such that $m\ne 0$
and $m(x+y)=m(x)\cd m(y)$ for every $x,y\in E$.
Since $m\in \Om$, we have
\begin{equation}\label{e1}
\sumjn c _j \cd  m(u\cd a_j )=0 \qquad (u\in S^1_a ).
\end{equation}
In the sequel we fix an exponential function $m$ with the properties above,
and look for a contradiction. We shall need the following result.

\begin{theorem}\label{tESS}
There exists a positive integer $A(N)$ that only depends on $N$ and has the
following property: whenever $n\le N$ and $\Ga$ is a multiplicative subgroup
of $\cc ^*$ of rank $r\le N$, then the number of solutions of the equation
\begin{equation*}
x_1 +\ldots +x_n =1
\end{equation*}
such that $x_1 \stb x_n \in \Ga$ and no subsum of $x_1 +\ldots +x_n$ equals
zero is at most $A(N)$.
\end{theorem}
This is obtained by applying \cite[Theorem 1.1]{ESS} with $K=\cc$ and
$a_1 =\ldots =a_n =1$. According to the theorem cited, we may take
$A(N)=\text{exp} (6N^{3N} (N+1))$.
(See also \cite[Theorem 6.1.3]{EG}.)

\begin{lemma} \label{l1}
Let $n\ge 2$. Then there are positive integers $d$ and $D$ that only depend on
$n$ and have the following property. Whenever $a_1 \stb a_n$, $c_1 \stb c_n$ and
$m$ are as above and satisfy \eqref{e1}, and $u\in S^1_a$, then there are
indices $1\le j_1 ,j_2 \le n$ such that
\begin{enumerate}[$\bullet$]
\item $m(u\cd (a_{j_2}-a_{j_1}) /d)$ and $m(u\cd i\cd (a_{j_2}-a_{j_1}) /d)$ are
roots of unity of degree dividing $D$, and
\item at least one of $m(u\cd (a_{j_2}-a_{j_1}) /d)$ and $m(u\cd i\cd
(a_{j_2}-a_{j_1}) /d)$ is different from $1$.
\end{enumerate}
\end{lemma}

\proof
It is enough to prove the statement for $u=1$. Indeed, suppose this special
case is true, and let $u_0 \in S^1_a$ be arbitrary. Replace $a_1 \stb a_n$
by $u_0 a_1 \stb u_0 a_n$. Note that the field generated by $u_0 a_1 \stb
u_0 a_n$ and $S^1_a$ is the same $E$ as above. Since $u_0 \cd u\in S^1_a$ for
every $u\in S^1_a$, it follows that $u_0 a_1 \stb u_0 a_n$, $c_1 \stb c_n$ and
$m$ satisfy \eqref{e1} with every $u\in S^1_a$. Therefore,
applying the special case $u=1$ for this system, we obtain the statement
for $u_0$.

Now we prove the statement with $u=1$.
We put $\ga _k =((1-k^2 ) +i\cd 2k )/(1+k^2 )$ for every $k=1,2,\ldots$.
Then $\ga _k \in S^1_a$ for every $k$.

For every $k$ there exists a partition $\{ 1\stb n\} =I_1 \cup \ldots \cup I_p$
with the following property: for every $1\le \mu \le p$,
\begin{equation*}
\sum_{j\in I_\mu} c_j \cd  m(\ga _k \cd a_j ) =0,
\end{equation*}
and $\sum_{j\in I} c_j \cd  m(\ga _k \cd a_j ) \ne 0$ whenever $\emp \ne I \subsetneq
I_\mu$. For a given $k$ there can be more than one such partition; we select one
for each $k$, and denote it by $\ii _{\ga _k}$.

Let $P(n)$ denote the number of partitions of $\{ 1\stb n\}$, and put
$B(n)=2\cd P(n) \cd A(3n) +1$. Then there is a set $H\su \{ 1\stb B(n)\}$
such that $|H|>2\cd A(3n)$, and the partitions $\ii _{\ga _k}$ $(k\in H)$ are the
same. Let $\ii _{\ga _k} =\{ I_1 \stb I_p\}$ for every $k\in H$.

Let $d= (1+B(n)^2 )!$. Then $d$ is a common multiple of the numbers $1+k^2$
$(k\in H )$, and thus $\ga _k =(e_k +i\cd f_k )/d$ for every $k\in H$, where
$|e_k | ,|f_k |\le d$. Note that $e_k \ne e_\ell$, if $k\ne \ell$.
Let $\mu \in \{ 1\stb p\}$ be given. Then, for every $k\in H$ we have
\begin{equation}\label{e3}
\begin{split}
0&=\sum_{j\in I_\mu} c_j \cd  m(\ga _k \cd a_j )=\sum_{j\in I_\mu} c_j \cd
m\left( e_k \cd \frac{a_j}{d} +f_k \cd \frac{i\cd a_j}{d} \right)\\
&=\sum_{j\in I_\mu} c_j \cd m(a_j /d)^{e_k} \cd  m(i\cd a_j /d)^{f_k} =
\sum_{j\in I_\mu} c_j \cd  u_j^{e_k} \cd v_j^{f_k} ,
\end{split}
\end{equation}
where $u_j =m(a_j /d)$ and $v_j =m(i\cd a_j /d)$. Note that $u_j , v_j \ne 0$
for every $j$ by the properties of the exponential.
Select an index $j_\mu \in I_\mu$. Then, by \eqref{e3}, we have
\begin{equation}\label{e4}
\sum_{j\in I_\mu , \ j\ne j_\mu}  \be _j \cd  (u_j /u_{j_\mu} )^{e_k} \cd
(v_j /v_{j_\mu} )^{f_k} =1
\end{equation}
for every $k\in H$, where $\be _j =-c_j /c _{j_\mu}$. (Recall that $c_j \ne 0$
for every $j=1\stb n$ by assumption.) Put $\ol u _j =u_j /u_{j_\mu}$ and
$\ol v _j =v_j /v_{j_\mu}$ $(j\in I_\mu )$, and let $\Ga$ be the multiplicative
group generated by the elements $\be _j , \ol u _j$ and $\ol v _j$. Then the
rank of $\Ga$ is at most $3n$, and $\be _j \cd   \ol u _j^{e_k} \cd \ol v _j^{f_k}
\in \Ga$ for every $j\in I_\mu$ and $k\in H$. By the choice of $A(3n)$,
the equation
\begin{equation}\label{e5}
\sum_{j\in I_\mu , \ j\ne j_\mu} x_j =1
\end{equation}
has at most $A(3n)$ solutions having the property that $x_j \in \Ga$ for every
$j$, and no subsum of the left hand side of \eqref{e5} is zero. However,
\eqref{e4} gives such a solution for every $k\in H$ by the choice of
$I_{\ga _k}$. Since $|H |> 2\cd A(3n)$, there must exist three distinct
indices $s,t,z\in H$ giving the same solution. Then
\begin{equation*}
\ol u _j ^{e_s} \cd \ol v _j^{f_s} = \ol u _j ^{e_t} \cd \ol v _j^{f_t} =
\ol u _j ^{e_z} \cd \ol v _j^{f_z}
\end{equation*}
for every $j\in I_\mu$, $j\ne j_\mu$. The equations above are also true if
$j= j_\mu$, since $\ol u _{j_\mu}= \ol v _{j_\mu} =1$. Then we have
\begin{equation}\label{e6}
  \ol u _j ^{e_t -e_s} \cd \ol v _j^{f_t -f_s} =1 \ \text{and} \ \ \ol u _j ^{e_z -e_s}
\cd  \ol v _j^{f_z -f_s} =1 \qquad (j\in I_\mu ).
\end{equation}
From \eqref{e6} we obtain $\ol u _j ^C =1$ and $\ol v _j ^C =1$ for every
$j\in I_\mu$, where
$$C=(e_z -e_s )(f_t -f_s ) -(e_t -e_s )(f_z -f_s ) .$$
We show that $C\ne 0$.
We have $e_k +i\cd f_k =d\cd \ga _k$, and $\ga _k \in S^1$ for every $k$.
Therefore, the points $(e_s ,f_s ), (e_t ,f_t ), (e_z ,f_z )$ are distinct,
and lie on a circle of radius $d$. Consequently, they are not collinear;
that is,
$$\frac{f_t -f_s}{e_t -e_s} \ne \frac{f_z -f_s}{e_z -e_s}.$$
Multiplying by the denominators we obtain $C\ne 0$. Note that
$|C|\le 8\cd d^2$.

We find that $\ol u _j$ and $\ol v _j$ are roots of unity of order at most
$|C| \le 8\cd d^2$ for every $j\in I_\mu$ and for every $\mu =1\stb p$.
Putting $D=(8\cd d^2 )!$, the orders
of $\ol u _j$ and $\ol v _j$ will be divisors of $D$.

Now we prove that there is an index $j\in \{ 1\stb n\}$ such that at least
one of $\ol u _j$ and $\ol v _j$ is different from $1$. Suppose not. Then,
for every $\mu =1\stb p$ we have $\ol u _j =\ol v _j =1$ for every $j\in I_\mu$.
By \eqref{e4}, we have $\sum_{j\in I_\mu , \ j\ne j_\mu}  \be _j =1$
and $\sum_{j\in I_\mu}  c_j =0$ $(\mu =1\stb p)$. However, this would imply
$\sumjn c_j =\sum_{\mu =1}^p \sum_{j\in I_\mu}  c_j =0$, which is impossible.

Therefore, we can find a $\mu$ and a $j\in I_\mu$ such that $\ol u _j =
u_j /u_{j_\mu} \ne 1$ or $\ol v _j = v_j /v_{j_\mu} \ne 1$. Choosing
$j_1 =j_\mu$ and $j_2 = j$ we find that  $m((a_{j_2}-a_{j_1}) /d) =
\ol u _j \ne 1$ or $m(i\cd (a_{j_2}-a_{j_1}) /d) =\ol v _j \ne 1$,
completing the proof. \hfill $\square$

\begin{lemma}\label{l2}
Let $S^1_a =A_1 \cup \ldots \cup A_N$ be a cover of $S^1_a$, and let
$c>1$ be an integer. Then there is a $j\in\{1\stb N\}$ and there are
elements $u_1 ,u_2 ,u_3 \in A_j$ and integers $n_1 ,n_2 ,n_3$ such that
$n_1 u_1 +n_2 u_2 +n_3 u_3 =0$ and $n_1 +n_2 +n_3$ is prime to $c$.
\end{lemma}
\proof
The polynomial $p(x)=cx^2 +x+c$ is irreducible over $\qq$, and its roots belong
to $S^1_a$. Let $\al$ be one of the roots of $p(x)$. Since $\al ^n \in S^1 _a$
for every $n$, there is a $j\in\{1\stb N\}$ such that $\al ^n \in A_j$ holds
for at least three distinct nonnegative exponents $n$. Suppose $\al ^r ,
\al ^s , \al ^t \in A_j$, where $0\le r<s<t$ are integers.

Since $\al ^r , \al ^s , \al ^t \in \qq (\al )$ and $\qq (\al )$ is a linear
space of dimension two over $\qq$, there are rational numbers
$n_1 ,n_2 ,n_3$, not all zero, such that $n_1 \al ^r +n_2 \al ^s  +n_3 \al ^t =0$.
Then $\al$ is a root of the polynomial $n_1 x^r +n_2 x^s  +n_3 x^t$, hence we
have
\begin{equation}\label{e7}
n_1 x^r +n_2 x^s  +n_3 x^t =(cx^2 +x+c)\cd q(x),
\end{equation}
where $q$ is a polynomial with rational coefficients. Let $q(x)=\sum_{i=u}^v
b_i x^i$, where $u\le v$ and $b_u \ne 0$, $b_v \ne 0$. Multiplying by the
common denominator of the coefficients $b_i$, we may assume that
$b_u \stb b_v$ are integers, and the polynomial $q$ is primitive, meaning
that the greatest common divisor of $b_u \stb b_v$ is $1$. Then
$n_1 ,n_2 ,n_3$ are integers. Since $cx^2 +x+c$ is also primitive, it follows
from Gauss' lemma that $n_1 x^r +n_2 x^s  +n_3 x^t$ is primitive as well.
It follows from \eqref{e7} that either $n_3 =0$ or $n_3 =c\cd b_v$.
In both cases we have $c\mid n_3$. We obtain $c\mid n_1$ similarly. Then $n_2$
must be prime to $c$, and thus the same is true for
$n_1 +n_2 +n_3$. \hfill $\square$

\subsection*{Conclusion of the proof of Theorem \ref{t1a}} Let $d$ and $D$ be as
in Lemma \ref{l1}. By Lemma \ref{l1}, for every $u\in S^1_a$ there are
indices $j_1 ,j_2$ such that $m(u\cd (a_{j_2}-a_{j_1}) /d)$ and
$m(u\cd i\cd (a_{j_2}-a_{j_1}) /d)$ are roots of unity of degree dividing $D$,
and at least one of them is different from $1$; that is, at least one
of them is the form $e^{2\pi i k /D}$ for some $0<k<D$. Therefore, we have
\begin{equation}\label{e9}
S^1_a = \bigcup_{j_1 =1}^n \bigcup_{j_2 =1}^n \bigcup_{k=1}^{D-1} (A_{j_1 , j_2 ,k}
\cup B_{j_1 ,j_2 ,k} ),
\end{equation}
where
$$A_{j_1 , j_2 ,k} =\{ u\in S^1_a \colon m(u\cd (a_{j_2} - a_{j_1}) /d ) =
e^{2\pi i\cd k/D}\} $$
and
$$B_{j_1 ,j_2 ,k} =\{ u\in S^1_a \colon m(u\cd i\cd (a_{j_2} - a_{j_1}) /d ) =
e^{2\pi i\cd k/D}\} .$$
Note that in \eqref{e9}, the index $k$ runs from $1$ to $D-1$, hence
$D\nmid k$.
Then, by Lemma \ref{l2}, there are elements $u_1 ,u_2 ,u_3$ and integers
$n_1 ,n_2 ,n_3$ such that $n_1 u_1 +n_2 u_2 +n_3 u_3 =0$, $n_1 +n_2 +n_3$
is prime to $D$, and $u_1 ,u_2 ,u_3$ belong to one of the sets
$A_{j_1 ,j_2 ,k}$ and $B_{j_1 , j_2 ,k}$.

Suppose they belong to $A_{j_1 ,j_2 ,k}$. We have $\sum_{t =1}^3 n_t \cd u_t
\cd (a_{j_2} - a_{j_1}) /d=0$, and thus
\begin{align*}
1&=m\left( \sum_{t =1}^3 n_t \cd u_t \cd (a_{j_2} - a_{j_1}) /d \right) =
\prod_{t =1}^3 m(u_t \cd (a_{j_2} - a_{j_1}) /d )^{n_t} = \\
&=\left( e^{2\pi i\cd k/D} \right) ^{n_1 +n_2 +n_3} =e^{2\pi i\cd k\cd (n_1 +n_2 +n_3 )/D} .
\end{align*}
This implies $D\mid k\cd (n_1 +n_2 +n_3 )$. However, $n_1 +n_2 +n_3$ is prime
to $D$ and $D\nmid k$, which is a contradiction. If
$u_1 ,u_2 ,u_3 \in B_{j_1 ,j_2 ,k}$, then we reach a contradiction by a similar
computation. \hfill $\square$

\section{Proof of Theorem \ref{t1a} for $k>2$}
We prove the statement by induction on $k$. By the results of the previous
section, the statement of the theorem is true for $k=2$.
Let $k\ge 2$, and suppose that the statement is true in $\rkk$. We prove
the statement in $\rr^{k+1}$. Let $a_1 \stb a_n \in \rk1$ be given.
Since the statement is obvious if $a_1 =\ldots =a_n$,
we may assume that $n\ge 2$ and $a_1 \ne a_n$. Using a translation
we may also assume that $a_1 =0$ and $a_n \ne 0$.

Let $S^k$ denote the unit sphere in $\rk1$; that is, let $S^k =\{ x\in
\rk1 \colon |x|=1\}$. If $v\in S^k$, then we denote by $v^\perp$ the linear
subspace of $\rk1$ of dimension $k$ and perpendicular to $v$. If $V$ is a
linear subspace of $\rk1$, then we denote by $G(V)$ the family of rigid motions
of $V$ as a Euclidean space. More precisely, if $d$ is the dimension of $V$,
then $G(V)$ is the family of
maps $\pi ^{-1} \circ \phi \circ \pi$, where $\pi$ is a fixed isometry mapping
$V$ onto $\rr ^d$, and $\phi \in G(\rr ^d )$. Thus $G_{k+1} =G(\rk1 )$.

\subsection*{The idea of the proof} Let $c_1 \stb c_n \in \cc$ be given
such that $\sumjn c_j \ne 0$. Suppose we only want to exclude the existence of
an exponential function $m\colon \rk1 \to \cc$ satisfying $\sum_{j=1}^n c_j
\cd m(\phi (a_j )) =0$ for every $\phi \in G(\rk1 )$. In order to prove this
it is enough to find a unit vector $v$ such that $\sum_{j=1}^n c_j \cd
m(t_j v )\ne 0$, where $t_j =\lan v,a_j \ran$ $(j=1\stb n)$. Indeed, suppose
we have found such a vector $v$. Every element of $\rk1$ has a unique
representation of the form
$b+tv$, where $b\in v^\perp$ and $t\in \rr$. Let $a_j =b_j +t_j v$ $(j=1\stb n)$.
If $\psi \in G(v^\perp )$, then the map $\ol \psi (b+tv)=\psi (b)+tv$ is a
rigid motion of $\rk1$, and thus
$$0=\sumjn c_j \cd m(\ol \psi (a_j ))= \sumjn c_j \cd m(\psi (b_j ) +t_j v)=
\sumjn c_j \cd m(\psi (b_j ))\cd m(t_j v).$$
Putting $d_j = c_j \cd m(t_j v )$ $(j=1\stb n)$, this implies $\sumjn d_j \cd
m(\psi (b_j ))=0$ for every $\psi \in G(v^\perp )$. Since $\sumjn d_j \ne 0$,
the induction hypothesis gives $m=0$ on $v^\perp$, which is impossible, as $m$
is nowhere zero.

Now we look for a vector $v$ with the desired properties, maybe not for the
set $A=\{ a_1 \stb a_n \}$, but for a congruent copy of $A$. Let $v_0$ be an
arbitrary unit vector, let $t_j =\lan v_0 ,a_j \ran$ $(j=1\stb n)$, and put
$B=\{ t_j v_0 \colon j=1\stb n\}$. If $V_0$ is a $k$-dimensional subspace of
$\rk1$
containing $v_0$ then, by $m\ne 0$ and by the induction hypothesis, we obtain
a rigid motion $\phi \in G(V_0 )$ such that $\sumjn c_j \cd m(\phi (t_j v_0 ))
\ne 0$. Using the multiplicative property of $m$, we can see that $\phi$ can
be chosen in such a way that $\phi (0)=0$ holds, and then $\phi$ is a linear
transformation preserving scalar products. Putting $v=\phi (v_0 )$ this implies
$\lan v, \phi (a_j )\ran = \lan v_0 ,a_j \ran =t_j$ for every $j$. Then the
argument above,  with $\phi (A)$ in place of $A$, leads to a contradiction.

This argument does not prove Theorem \ref{t1}, since harmonic analysis does not
hold on the discrete additive group of $\rk1$, so the existence of a
counterexample to the statement of the theorem does not imply the existence
of a counterexample which is an exponential. However, as harmonic analysis
holds on countable Abelian groups, we can find a suitable countable subgroup
of $\rk1$ on which the previous argument can be implemented.

\medskip
Now we turn to the details of the proof of Theorem \ref{t1a}.
Let a unit vector $v_0 \in S^k$ be selected such that
$\lan v_0 ,a_n \ran \ne 0$, and put $t_j =\lan v_0 ,a_j \ran$ $(j=1\stb n )$.
Note that $t_1 =\lan v_0 ,0\ran =0$ and $t_n \ne 0$. Let $V_0$ be a
linear subspace of $\rk1$ of dimension $k$ containing $v_0$.
By the induction hypothesis applied to the points $t_1 v_0 \stb t_n v_0
\in V_0$, we find a countable additive group $E_0 \su V_0$ containing
$t_1 v_0 \stb t_n v_0$ and
having the following property: whenever $c_1 \stb c_n \in \cc$, $\sumjn c_j
\ne 0$, and $f\colon E_0 \to \cc$ is such that  $\sumjn c_j \cd
f(\phi (t_j v_0 ))=0$ for every $\phi \in G(V_0 )$ satisfying $\phi (t_j v_0 )
\in E_0$ for every $j=1\stb n$, then $f\equiv 0$.

Let $W=\{ v\in S^k \colon t_j v\in E_0 \ (j=1\stb n)\}$. Since $t_n \ne 0$,
$W$ is a countable set of unit vectors. For every $v\in W$ let a rigid motion
$\phi _v \in G_{k+1}$
be selected such that
\begin{equation}\label{eph}
  \phi _v (0)=0 \ \text{ and} \ v=\phi _v (v_0 ).
\end{equation}
Then $\phi _v$ is a linear transformation of $\rk1$. Let $b_{v,j}$ denote the
orthogonal projection of $\phi _v (a_j )$ onto $v^\perp$ $(j=1\stb n)$.

Let $v\in W$.
Applying the induction hypothesis again, we find a countable additive group
$E_v \su v^\perp$ containing $b_{v,1} \stb b_{v,n}$ and having the following
property: whenever $d_1 \stb d_n \in \cc$,
$\sumjn d_j \ne 0$, and $f\colon E_v \to \cc$ is such that  $\sumjn d_j \cd
f(\psi (b_{v,j} ))=0$ for every $\psi \in G(v^\perp )$ satisfying $\psi (b_{v,j} )
\in E_v$ for every $j=1\stb n$, then $f(b_{v,1} )=0$.

Let $E$ be the additive group generated by $E_0 \cup \bigcup_{v\in W}
\bigcup_{j=1}^n (E_v +t_j v)$. Then $E$ is countable. We show that $E$ satisfies
the condition of Theorem \ref{t1a}.

Let $c_1 \stb c_n$ be fixed complex numbers satisfying $\sumjn c_j \ne 0$,
and let $\Lambda$ denote the set of all functions $f\colon E \to \cc$ such that
$\sumjn c_j \cd f(\phi (a_j ))=0$ for every $\phi \in G_{k+1}$ satisfying
$\phi (a_j )\in E$ for every $j=1\stb n$. It is clear that $\Lambda$ is a linear
space over $\cc$, and that $\Lambda$ is invariant under translations
by elements of $E$. It is also easy to see that $\Lambda$ is closed in the set
$\cc ^E$ equipped with the product topology. This means that $\Lambda$ is a
variety on the discrete additive group $E$.

Suppose  $\Lambda\ne \{ 0\}$. Then, by Theorem \ref{tha}, $\Lambda$ contains
an exponential; that is, a function $m\colon E\to \cc$ such that $m\ne 0$
and $m(x+y)=m(x)\cd m(y)$ for every $x,y\in E$. Since $m\in \Lambda$, we have
$\sumjn c_j \cd  m(\phi (a_j ))=0$ for every $\phi \in G_{k+1}$ such that
$\phi (a_j )\in E$ $(j=1\stb n)$.

Now $m$ is defined on $E_0$, and $m$ is nowhere zero. Then it follows from the
choice of $E_0$ that there exists a $\phi \in G(V_0 )$ satisfying
$\phi (t_j v_0 )\in E_0$ for every $j=1\stb n$, and such that
$\sumjn c_j \cd m(\phi (t_j v_0 ))\ne 0$. Since $t_1 =0$, we have
$\phi (0)\in E_0$. Put $\si =\phi -\phi (0)$. Then $\si \in G(V_0 )$
and $\si (t_j v_0 )\in E_0$ for every $j=1\stb n$, as $E_0$ is an additive
group. Note that $\si$ is a linear transformation of $V_0$.
We put $d_j =c_j \cd m(\si (t_j v_0 ))$ $(j=1\stb n)$. Then $\sumjn d_j \ne 0$,
since $m(\si (t_j v_0 )) =m(\phi (t_j v_0 ))/ m(\phi (0))$ for every $j$.

Let $v=\si (v_0 )$. Then $v\in W$, as $t_j  v=t_j \si (v_0 )=
\si (t_j v_0 )\in E_0$ for every $j=1\stb n$.
We show that if $\psi \in G(v^\perp )$ is such that $\psi (b_{v,j} ) \in E_v$
for every $j=1\stb n$, then
\begin{equation}\label{evege}
\sumjn d_j \cd m(\psi (b_{v,j} ))=0.
\end{equation}
As $m$ is defined on $E_v$, and is nowhere zero, this will contradict the
choice of $E_v$, proving the theorem.

Every element of $\rk1$ has a unique representation of the form $b+tv$,
where $b\in v^\perp$ and $t\in \rr$. Putting $\ol \psi (b+tv )=\psi (b)+tv$,
we define the rigid motion $\ol \psi \in G_{k+1}$. We prove that
$\phi _v (a_j ) =b_{v,j} +t_j v$ for every $j=1\stb n$. (As for $\phi _v$,
see \eqref{eph}.) Indeed, if $\phi _v (a_j ) =b_{v,j} +tv$, then
$$t=\lan v,\phi _v (a_j ) \ran =\lan \phi _v (v_0 ),\phi _v (a_j ) \ran =
\lan v_0 , a_j \ran =t_j .$$
Therefore, we have $(\ol \psi \circ \phi _v )(a_j )= \psi (b_{v,j} )+
t_j v \in E_v +t_j v \su E$, and
$$m((\ol \psi \circ \phi _v )(a_j ))= m(\psi (b_{v,j} ))\cd m(t_j v).$$
Now $d_j =c_j \cd m(\si (t_j v_0 ))=c_j \cd m(t_j v)$, and thus
$$\sumjn d_j \cd m(\psi (b_{v,j} )) =\sumjn c_j \cd m(t_j v) \cd
m(\psi (b_{v,j} )) =   \sumjn c_j \cd m((\ol \psi \circ \phi _v )(a_j ))=0,$$
as $\ol \psi \circ \phi _v \in G_{k+1}$, and $(\ol \psi \circ \phi _v )
(a_j )\in E$ for every $j=1\stb n$. This completes the proof
of \eqref{evege} and of the theorem. \hfill $\square$

\section{On infinite systems of linear equations}
Let $K$ be a field, let $I$ be a nonempty set, and let $\ff$ denote
the set of functions $f\colon I\to K$ such that $\{ i\in I\colon f(i)\ne 0\}$
is finite. By a {\it system of linear equations over $K$} we mean a subset $\ee$
of $\ff \times K$. A {\it solution of the system $\ee$} is a function
$x\colon I\to K$ such that $\sum_{i\in I} f(i)\cd x(i)=b$ for every
$(f,b)\in \ff \times K$.
\begin{proposition}
If every finite subsystem of $\ee$ has a solution, then so
has $\ee$.
\end{proposition}
\proof
Let $V$ denote the set of linear combinations of the elements
of $\ee$; then $V$ is a linear subspace of $\ff \times K$ as a vector space
over $K$. We prove that if $(0,b)\in V$, then $b=0$. Indeed, let $(0,b)=
\sum_{j=1}^n c_j \cd (f_j ,b_j )$, where $(f_j ,b_j )\in \ee$ and $c_j \in K$
$(j=1\stb n)$. By assumption, the finite system $\{ (f_j ,b_j )\colon
j=1\stb n\}$ has a solution $x$. Then $\sum_{i\in I} f_j (i)\cd x(i)=b_j$ for
every $j=1\stb n$. Taking the linear combination of these equations with
coefficients $c_j$, we obtain $0=b$, as we stated.

This implies that if $(f,b)$ and $(f,c)$ are both elements of $V$, then $b=c$.
Putting $L(f)=b$ for every $(f,b)\in V$ we define a linear map on
the linear subspace $\{ f\colon (\exists b) ((f,b)\in V)\}$ of $\ff$.
Since $L$ is linear, it can be extended to a
linear map $\ol L \colon \ff \to K$. In particular, $\ol L$ is defined
on the characteristic function $\chi _{\{ i\}}$ of the singleton $\{ i\}$
for every $i\in I$. It is easy to check that $x(i)=\ol L
(\chi _{\{ i\}} )$ $(i\in I)$ defines a solution to the system
$\ee$. \hfill $\square$

\vspace*{14pt}

\section*{Acknowledgment}
We are grateful to the anonymous reviewer for her/his detailed comments and
suggestions, which significantly improved the quality of the article.
We also thank Andr\'as M\'ath\'e for his valuable comments on an earlier
version of the paper.

\paragraph{Funding statement}
Both authors were supported by the Hungarian National Foundation for
Scientific Research, Grant No. K146922. The first author
was also supported by the J\'anos Bolyai Research Fellowship,
ÚNKP-23-5-ELTE-1275
New National Excellence Program of the Ministry for Culture and
Innovation from the source of the National Research, Development and Innovation Fund, OTKA grant no. FK 142993 and Starting 150576.

\begin{small}\noindent
(G. Kiss)\\
{\sc Corvinus University of Budapest\\
F\H ov\'am t\'er 13-15, Budapest, H-1093, Hungary\\
and\\
HUN-REN Alfred Renyi Mathematical Institute\\
Reáltanoda u. 13-15, Budapest, H-1053, Hungary\\
E-mail: {\tt kiss.gergely@renyi.hu, kigergo57@gmail.com}}

(M. Laczkovich)\\
{\sc Department of Analysis, E\"otv\"os Lor\'and University\\
P\'azm\'any P\'eter s. 1/C, Budapest, H 1117, Hungary\\
E-mail: {\tt miklos.laczkovich@gmail.com}}
\end{small}


\begin{thebibliography}{99}

\bibitem{ESS} J.-H. Evertse, H.P. Schlickewei and W.M. Schmidt, Linear
equations in variables which lie in a multiplicative group, {\it Ann. of Math.}
(2) {\bf 155} (2002), no. 3, 807-836.

\bibitem{EG} J.-H. Evertse and K. Győry: {\it Unit equations in Diophantine
number theory.} Cambridge Stud. Adv. Math., 146. Cambridge University Press,
Cambridge, 2015.

\bibitem{GMW} S. Gao, A.W. Miller and W.A.R. Weiss,
Steinhaus sets and Jackson sets {\it Advances in logic}, 127-145.
Contemp. Math., 425.

\bibitem{GD} C. de Groote, M. Duerinckx, Functions with constant mean on similar countable subsets of $\sik$, {\it Amer. Math. Monthly}  {\bf 119} (2012), 603-605.

\bibitem{HJL} D. Henkis, S. Jackson, J. Lobe, The Finite Steinhaus Problem,
{\it The Quarterly Journal of Mathematics} {\bf 67}(4) (2016), 551-564.

\bibitem{JM1} S. Jackson and R.D. Mauldin, On a lattice problem of H. Steinhaus
{\it J. Amer. Math. Soc.} {\bf 15} (2002), no.4, 817-856.

\bibitem{JM} S. Jackson and R.D. Mauldin, Survey of the Steinhaus tiling
problem, {\it Bull. Symbolic Logic} {\bf 9} (2003) no.3, 335-361.

\bibitem{KKS} R. Katz, M. Krebs and A. Shaheen, Zero sums on unit square vertex
sets and plane colorings, {\it Amer. Math. Monthly} {\bf 121} (2014), no. 7,
610-618.

\bibitem{KLV} G. Kiss, M. Laczkovich and Cs. Vincze, The discrete Pompeiu
problem on the plane, {\it Monatsh. Math.} {\bf 186} (2018), no. 2, 299-314.

\bibitem{KMS} G. Kiss, R. D. Malikiosis, G. Somlai and M. Vizer, On the
discrete Fuglede and Pompeiu problems, \textit{Analysis \& PDE}, \textbf{13}(3)
(2020), 765-788.

\bibitem{KP} M. Kolountzakis and M. Papadimitrakis, Measurable Steinhaus sets
do not exist for finite sets or the integers in the plane,  {\it Bulletin of
the London Mathematical Society} {\bf 49}(5) (2017), 798-805.

\bibitem{KS} M. Kolountzakis and E. Spyridakis, Curves in the Fourier zeros
of polytopal regions and the Pompeiu problem, {\it Anal. Math.} {\bf 50}
(2024), no. 4, 1081-1098.

\bibitem{K1} P. Komj\'ath, A lattice-point problem of Steinhaus, {\it Quart.
J. Math. Oxford} Ser. (2) 43 (1992), no. 170, 235-241.

\bibitem{K2} P. Komj\'ath, A coloring result for the plane, {\it J. Appl.
Anal.} {\bf 5} (1999), no. 1, 113-117.

\bibitem{LS} M. Laczkovich and G. Sz\'ekelyhidi, Harmonic analysis on
discrete Abelian groups, {\it Proc. Amer. Math. Soc.} {\bf 133}  (2005),
no. 6, 1581-1586.

\bibitem{LP} P. A. Linnell and M. J. Puls, The two-sided Pompeiu problem for
discrete groups, {\it Proc. Amer. Math. Soc.}, Series B, {\bf 9} (2-22)
(2022), 221-229.

\bibitem{MR} F.C. Machado and S. Robins, The null set of a polytope, and
the Pompeiu property for polytopes {\it J. Anal. Math.} {\bf 150} (2023), no.2,
673-683.

\bibitem{Puls} M. J. Puls, The Pompeiu problem and discrete groups,
{\it Monatsh. Math.} {\bf 172}(3-4) (2013), 415-429.

\bibitem{R} A. G. Ramm, The Pompeiu problem, {\it Applicable Analysis} {\bf 64}
no. 1-2 (1997), 19-26.

\bibitem{Sch} J.H. Schmerl, Coloring $\rr ^n$, {\it Trans. Amer. Math. Soc.}
{\bf 354} (2002), no. 3, 967-974.

\bibitem{S} L.E. Shader, All right triangles are Ramsey in $E^2$! {\it J.
Combinatorial Theory} Ser. {\bf A20} (1976), no.3, 385-389.

\bibitem{Za} L. Zalcman, A bibliographical survey of the Pompeiu Problem, in
the book {\it Approximation by solutions of partial differential equations.}
Edited by B. Fuglede. Kluwer Acad., Dordrecht, 1992, 177-186.

\bibitem{Z} D. Zeilberger, Pompeiu's problem on discrete space, {\it Proc.
Natl. Acad. Sci. USA}, {\bf 75}(8) (1978), 3555-3556.





\end{thebibliography}
\end{document}